\documentclass{elsart}
\usepackage{amsfonts}
\usepackage{amsmath}
\usepackage{amssymb}
\usepackage{amscd}
\usepackage{verbatim}
\usepackage{graphicx}
\usepackage{subfigure}
\usepackage[T1]{fontenc}
\usepackage[latin1]{inputenc}
\usepackage[english]{babel}

\input xy
\xyoption{all}

\def\path{\mathop{\rightsquigarrow}\nolimits}
\def\mod{\mathop{\rm mod}\nolimits}

\newcommand{\mor}[3]{\xymatrix@1@C=15pt{#1: #2\ar[r]& #3}}
\newcommand{\map}[2]{\xymatrix@1@C=15pt{#1\ar@{|->}[r]& #2}}

\begin{document}

\begin{frontmatter}

\title{Laura string algebras}%

\author[JulieAd]{Julie Dionne\thanksref{JulieISM}},
\thanks[JulieISM]{Partially supported by Université de Sherbrooke.}
\ead{julie.dionne2@usherbrooke.ca}
\address[JulieAd]{D\'epartement de math\'ematiques, Université de
Sherbrooke, 2500, boul. de l'Université, Sherbrooke, Qu\'ebec,
Canada, J1K 2R1}

\begin{abstract}
\noindent We give a simple combinatorial criterion allowing to recognize whether a string (or, more generally, a special biserial) algebra is a laura algebra or not. We also show that a special biserial algebra is laura if and only if it has a finite number of isomorphism classes of indecomposable modules which have projective dimension and injective dimension greater than or equal to two, solving a conjecture ok Skowro{\'n}ski for special biserial algebras.
\end{abstract}


\end{frontmatter}

\bibliographystyle{plain}

\maketitle

%-------------------------------------------------------------------------------
%
%  SECTION  : INTRODUCTION
%
%-------------------------------------------------------------------------------
%\section{Introduction}
%    \label{Introduction}

Let $k$ be an algebraically closed field. A finite-dimensional $k$-algebra $R$ is biserial if the radical of every projective indecomposable $R$-module  is the sum of two uniserial modules whose intersection is simple or zero \cite{F79}. In 1983, Skowro{\'n}ski and Waschbüsch characterized biserial algebras of finite representation type by the fact that almost split sequences have at most two non-projective middle terms \cite{SW83}. The biserial algebras which have at most two middle terms in their almost split sequences are called string algebras (see \cite{BR87}). The following definition is equivalent:
 \begin{defn}
\label{stringAlgebra}
A $k$-algebra $R$ is a string algebra if it admits a presentation $R=kQ/I$ such that : \begin{enumerate}
\item Each point has at most two arrows entering and two arrows exiting;
\item For each arrow  $\alpha:x \rightarrow y$ there is at most one arrow $\beta:y \rightarrow z $ such that $\alpha\beta$ is not in $I$ and at most one arrow $\gamma: z \rightarrow x$ such that $\gamma\alpha$ is not in $I$;
\item The ideal $I$ is monomial.
\end{enumerate}
\end{defn}
The characterizations of string algebras which are of finite representation type \cite{BR87}, tilted \cite{HL00}, quasi-tilted \cite{HL99} and shod \cite{BT05} are known. For a string algebra $R=kQ/I$, we have the following definitions. A walk in the quiver of $R$ is \emph{reduced} if it contains no subwalk of the form $\alpha\alpha^{-1}$ or of the form $\alpha^{-1}\alpha$. A \emph{cycle} is a non-oriented cycle, that is, a reduced walk starting and ending in the same point of the quiver. A walk $\omega$ is a \emph{string} if it is reduced and contains no relation. Given a string $\omega$, the representation admitting a copy of the field $k$ at $x$ for each passage of $\omega$ on $x$ and with the obvious morphisms is the \emph{string module} over $\omega$, denoted by $M(\omega)$.
A walk $\omega$ is a \emph{band} if it is a cyclic string which is not the power of another cyclic string and if there is no $n$ such that $\omega^n$ is in $I$.

Laura algebras have been defined independently by Assem and Coelho \cite{AC03} and Skowro{\'n}ski \cite{S03}.
Let $A$ be an algebra, its left part is the full subcategory of ind$A$ defined as follows:
$$\mathcal L_A=\{M \in \mbox{ind}A \textit{ }|\textit{ }\textit{ for all }\textit{ }
L \rightsquigarrow M,\textit{ }\mbox{dp}L \leq 1\}$$ Its right part is dually defined and is denoted by $\mathcal R_A$.

\begin{defn}
An algebra $A$ is laura if \emph{ind}$A \setminus (\mathcal L_A \cup
\mathcal R_A)$ contains finitely many objects.
\end{defn}

We say that $R$ is \emph{strict laura} if it is laura but not quasi-tilted.
Skowro{\'n}ski conjectured that an algebra $R$ is laura if and only if the number of indecomposable modules with projective and injective dimension greater than or equal to two is finite.

Our aim is to characterize laura string algebras and to show Skowro{\'n}ski's conjecture for these algebras.

\begin{defn}
Let $R \cong kQ/I$ be a string algebra, and
$\omega$ be a double-zero on $Q$ (see \cite{HL99} for instance). We say that $\omega$ is an interlaced double-zero, abbreviated by DOZE, if there exist a
band $\omega_2$, two walks $\omega_1$ and
$\omega_3$ and two relations $\rho_1$ and $\rho_2$ such that $\omega=\rho_1\omega_1{\omega_2}\omega_3\rho_2$.
\end{defn}

Remark that if $\omega=\rho_1\omega_1{\omega_2}\omega_3\rho_2$ is a DOZE, then $\rho_1\omega_1{\omega_2}^n\omega_3\rho_2$ is a double-zero for all $n \geq 0$.

The first two sections will be useful for proving the following theorem.
Recall that an algebra is quasi-tilted of canonical type if its bounded derived category is equivalent to the bounded derived category of a category of coherent sheaves coh$\mathbb X$ on a weighted projective line $\mathbb X$ in the sense of Geigle and Lenzing.

\begin{thm}\label{TheoA}
Let $R=kQ/I$ be a string algebra having no DOZE and such that $Q$ contains a band.
\begin{enumerate}
\item[a)] If $Q$ has at least one band with at least one arrow entering and at least one arrow exiting, then $R$ is quasi-tilted of canonical type.
\item[b)] If $Q$ contains at least one band and is such that each band has either only exiting arrows, or only entering arrows, then $R$ is strict laura or tilted.
\end{enumerate}
\end{thm}

In the third section, we prove the validity of Skowro{\'n}ski's conjecture for string and more generally for special biserial algebras. Our main results are the following two theorems.

\begin{thm}\label{TheoPresqueFinal} Let $R=kQ/I$ be a string algebra. The following are equivalent: \begin{enumerate} \item[a)] There is no DOZE on $(Q,I)$; \item[b)] $R=kQ/I$ is
laura; \item[c)] $R$ has only a finite number of isomorphism classes of indecomposable modules having projective and injective dimension greater than or equal to 2.
\end{enumerate}
\end{thm}

If $R$ is a special biserial algebra, let $J$ be the minimal ideal containing all paths lying in non-monomial relation. Then $R/J$ is a string algebra.

\begin{thm}\label{TheoFinal} Let $R=kQ/I$ be a special biserial algebra. The following are equivalent: \begin{enumerate} \item[a)] There is no DOZE on $(Q,I)$;   \item[b)] $R=kQ/I$ is
laura;   \item[c)] $R$ has only a finite number of indecomposable modules having a projective and an injective dimension greater or equal than 2; \item[d)] $R/J$ is
laura.
\end{enumerate}
\end{thm}
%-------------------------------------------------------------------------------
%
%  CHAPITRE: Preliminares
%
%-------------------------------------------------------------------------------
\section{Preliminaries}

A quiver $Q$ is a tuple ($Q_0$, $Q_1$, $s:Q_1\rightarrow Q_0$, $t:Q_1\rightarrow Q_0$). We call $Q_0$ the set of vertices, $Q_1$ the set of arrows, $s(\alpha)$ the source of the arrow $\alpha$ and $t(\alpha)$ the target of the arrow $\alpha$. The algebra $kQ$ is the $k$-vector space generated by all the paths on $Q$ with the multiplication defined by the composition of paths. If $I$ is an ideal of $kQ$, consider the algebra $kQ/I$. We have a complete set of idempotents of $kQ/I$ given by the trivial paths on each vertex $x$, denoted by $\varepsilon_x$.
For a finite dimensional, basic and connected $k$-algebra $R\cong kQ/I$, we denote by $\mod R$ the category of finite dimensional left $R$-modules. The indecomposable projective, injective and simple modules associated to $\varepsilon_x$ will be denoted respectively by $P_x$, $I_x$ and $S_x$.
A path in mod$R$ from $M$ to $N$ is a sequence
$$\xymatrix@R=10pt@C=10pt{
(*)M=M_0\ar[rr]^{f_1}&&M_1\ar[rr]&& \ldots&&\ar[rr]^{f_t}&&M_t=N}$$
of non-zero morphisms between indecomposable modules. The path $(*)$ is called sectional if, for any $i$, we have that $X_{i}\ncong\tau X_{i+2}$. A refinement of $(*)$ is a path
$$\xymatrix@R=10pt@C=10pt{
(*)M=X_0\ar[rr]^{g_1}&&X_1\ar[rr]&& \ldots&&\ar[rr]^{g_s}&&X_s=N}$$
with $s\geq r$ such that there is an order-preserving function $\sigma: \{1,$ $2,$ $3,$ $...,$ $t-1\} \rightarrow \{1,$ $2,$ $3,$ $...,$ $s-1\}$ such that $M_i \cong X_{\sigma(i)}$ for all $i$. The path $(*)$ is called sectionally refinable if it has a sectional refinement.
It is well known that an algebra is laura if and only if there is only a finite number of isomorphism classes of indecomposable modules on a path going from an injective module to a projective module.
Let $R$ be a strict laura or a tilted algebra then $R$ has a unique faithful nonsemiregular component which is quasidirected. Let $(_l {\Gamma_\lambda})_{\lambda \in \Lambda}$ be the left stable parts of the faithful non semiregular components (since $R$ is laura, $\Lambda$ is a finite set). Let $(_l {\Sigma_\lambda})_{\lambda \in \Lambda}$ be complete slices of each one of those stable parts (since each component has a finite number of orbits, such slices exist). For all $\lambda \in \Lambda$, we define $_\infty R_\lambda$ to be the full subcategory generated by the support of $ _l \Sigma_\lambda$. The left end algebra $_\infty R$ of $R$ is by definition the product of the $_\infty R_\lambda$. We define dually the right end algebra $R_\infty$ of $R$.

%-------------------------------------------------------------------------------
%
%  CHAPITRE: Algèbres articulées simplement connexes
%
%-------------------------------------------------------------------------------
\section{String algebras without DOZE}
    \label{Structure}
Throughout this section, let $R=kQ/I$ be a string algebra. The following lemma will be useful for the next two sections:

\begin{lem}\label{cyclesdisjoints}
If $R=kQ/I$ has no DOZE, then two bands of $R$ intersect in at most one point.
\end{lem}
Proof: We verify that for all cases with more than one point in common, we obtain a DOZE.
If $\omega_1$ and $\omega_2$ have a non trivial walk $u$ in common, we have :
$$
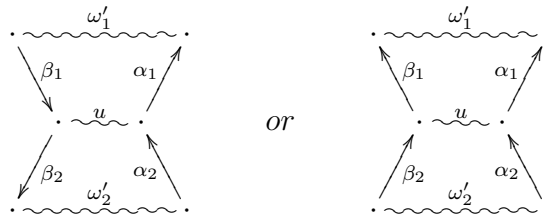
\xymatrix@R=10pt@C=10pt{
\cdot\ar[rdd]^{\beta_1} &&&& \cdot \ar@{~}[llll]_{\omega'_1}&&&&\cdot &&&& \cdot \ar@{~}[llll]_{\omega'_1}\\
&&&&&&&&&&&&\\
& \cdot\ar@{~}[rr]^u \ar[ldd]^{\beta_2}& &\cdot\ar[ruu]^{\alpha_1}&&&or&&  & \cdot\ar@{~}[rr]^u \ar[luu]_{\beta_1}&    &\cdot\ar[ruu]^{\alpha_1}&  \\
&&&&&&&&&&&&\\
\cdot &&&& \cdot \ar@{~}[llll]_{\omega'_2} \ar[luu]^{\alpha_2}
&&&&\cdot\ar[ruu]_{\beta_2} &&&& \cdot \ar@{~}[llll]_{\omega'_2}\ar[luu]^{\alpha_2}
}$$
where $\omega_1= u \alpha_1\omega'_1 \beta_1$ and $\omega_2=u {\alpha_2}^{-1}\omega'_2 {\beta_2}^{-1}$ in the first case and $\omega_1= u \alpha_1\omega'_1 {\beta_1}^{-1}$ and $\omega_2=u {\alpha_2}^{-1}\omega'_2 {\beta_2}$ in the second case.
In the first case, $\alpha_2\omega_1'\beta_2$ is a DOZE. In the second, $\alpha_2\omega_1' u \omega_2' \beta_1$ is a DOZE.
If $\omega_1$ and $\omega_2$ have no non trivial walks in common, we have one of the following three cases:
$$
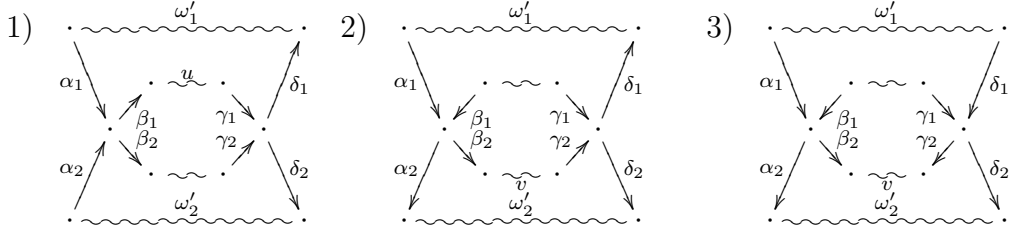
\xymatrix@R=8pt@C=8pt{
1)&\cdot\ar[rdd]_{\alpha_1} &&&&&& \cdot \ar@{~}[llllll]_{\omega'_1}&2)&\cdot\ar[rdd]_{\alpha_1} &&&&&& \cdot \ar@{~}[llllll]_{\omega'_1}&&3)&\cdot\ar[rdd]_{\alpha_1} &&&&&& \cdot \ar@{~}[llllll]_{\omega'_1}\ar[ldd]^{\delta_1}\\
&&&\cdot&&\cdot\ar@{~}[ll]_u\ar[rd]_{\gamma_1}&&&&&&\cdot\ar[dl]^{\beta_1}&&\cdot\ar@{~}[ll]\ar[rd]_{\gamma_1}&&
&&&&&\cdot\ar[dl]^{\beta_1}&&\cdot\ar@{~}[ll]\ar[rd]_{\gamma_1}&&\\
&& \cdot\ar[ru]_{\beta_1} \ar[rd]^{\beta_2}&&& &\cdot\ar[ruu]_{\delta_1}\ar[rdd]^{\delta_2} && && \cdot \ar[rd]^{\beta_2}\ar[ldd]_{\alpha_2}&&& &\cdot\ar[ruu]_{\delta_1}\ar[rdd]^{\delta_2} & &&&& \cdot \ar[rd]^{\beta_2}\ar[ldd]_{\alpha_2}&&& &\cdot\ar[rdd]^{\delta_2}\ar[dl]_{\gamma_2} & \\
&&&\cdot&&\cdot\ar@{~}[ll]\ar[ru]^{\gamma_2}&&&&&&\cdot&&\cdot\ar@{~}[ll]^v\ar[ru]^{\gamma_2}&&&&&&&\cdot&&\cdot\ar@{~}[ll]^v&&\\
&\cdot\ar[ruu]^{\alpha_2} &&&&&& \cdot \ar@{~}[llllll]_{\omega'_2}&&\cdot &&&&&& \cdot \ar@{~}[llllll]_{\omega'_2}& &&\cdot&&&&&& \cdot \ar@{~}[llllll]_{\omega'_2}
}$$
where $\omega_1= \beta_1 u \gamma_1 \delta_1\omega'_1 \alpha_1$ and $\omega_2=\beta_2 v \gamma_2 \delta_2\omega'_2 \alpha_2$ in the first case, $\omega_1= {\beta_1}^{-1} u \gamma_1 \delta_1\omega'_1 \alpha_1$ and $\omega_2=\beta_2 v \gamma_2 \delta_2\omega'_2 {\alpha_2}^{-1}$ in the second case and $\omega_1= {\beta_1}^{-1} u \gamma_1 {\delta_1}^{-1} \omega'_1 \alpha_1$ and $\omega_2=\beta_2 v {\gamma_2}^{-1} \delta_2\omega'_2 {\alpha_2}^{-1}$ in the last case.

Then, in the first case $\alpha_1\beta_1 u\gamma_1\delta_1$ is a DOZE,
in the second $\beta_1\beta_2 v\gamma_2\delta_2$ is a DOZE, and in the third $\beta_1\beta_2 v\gamma_2^{-1}\delta_1^{-1}w\alpha_1\alpha_2$ is a DOZE.
\\$\square$

\subsection{Quasi-tilted algebras without DOZE}
\begin{rem}
Let $R=kQ/I$ be a string algebra without DOZE having a band $\Theta$ with entering arrow $\beta$ and exiting arrow $\alpha$. Let $\beta^+$ be an arrow of $\Theta$ such that $\beta \beta^+ \in I$. Then there is a string $\omega$ on $\Theta$, going from the source of $\beta^+$ to the source of $\alpha$, such that $\omega \alpha$ is a string.
\end{rem}

\begin{prop}
Let $Q$ be a quiver with a band having at least one entering arrow $\beta$ and at least one exiting arrow
$\alpha$. If $R=kQ/I$ is a string algebra having no DOZE, then $R$ has no double-zero.
\end{prop}

Proof: We denote by $\beta^+$ the arrow such that $\beta\beta^+ \in I$ and by $\alpha^-$ the arrow such that $\alpha^- \alpha \in I$. Suppose that we have a double-zero. If $\beta\beta^+$ is the relation at the beginning of the double-zero then we have a DOZE. Thus, suppose that this is not the case. Then we can find a walk of minimal length between the end point of $\beta^+$ and a point on the double-zero. If the first relation on the walk of minimal length composed with a part of the double-zero points in the same direction than $\beta\beta^+$, we have a DOZE. Otherwise, let $\omega$ be the walk going from $t(\beta)$ to $s(\alpha^-)$ such that $\beta\omega$ is a string. Then we have a DOZE.
\\$\square$

\subsection{Cycles on a string algebra without DOZE}
We now consider the case where each cycle has only entering arrows or only exiting arrows. We ignore the trivial case where $R$ has only one band with neither entering nor exiting arrow (in this case, $R$ is hereditary). We denote by $\Theta_1, \Theta_2, ..., \Theta_n$ the bands (up to a cyclic permutation) having only exiting arrows. Since $R$ is finite dimensional, a band $\Theta$ has a point $a$ such that all arrows starting at $a$ belong to $\Theta$. We arbitrarily choose such a point on each cycle $\Theta_i$ and denote it by
$a_i$. Let $\alpha$ be an exiting arrow of $\Theta_i$. We denote by
$\alpha^-$ the arrow of $\Theta_i$ such that $\alpha^-\alpha \in I$. We denote by $\rho'_{(i,\alpha)}$
the minimal reduced walk going from $a_i$ to $t(\alpha)$ and passing through
$\alpha^-\alpha$.

We denote by $\Delta_1, \Delta_2, ..., \Delta_m$ the bands (up to a cyclic permutation) having only entering arrows.

\begin{defn} Let $R=kQ/I$ and $\omega$ be a string of $(Q, I)$. We denote by $W(\omega)$
the set of strings $\omega'$ of $(Q, I)$ such that $\omega'
= \omega_1\omega\omega_2$. We denote by $D(\omega)$ the subcategory of
$(Q, I)$ whose objects are the points $x$ such that there exists
a string in $W(\omega)$ passing through $x$ and whose morphisms are the composition of arrows $\alpha$ for which there exists a string of $W(\omega)$ passing through $\alpha$.
\end{defn}

\begin{defn}
Let $R=kQ/I$ be a string algebra such that each band on its quiver has only entering arrows or only exiting arrows. We define $A_i$ to be the subcategory of $R$ whose objects are the points of $D(\varepsilon_{a_i})$ and whose morphisms are given by the linear combinations of paths of arrows of $D(\varepsilon_{a_i})$. The algebra $A_i$ is the quotient of the path algebra of $D(\varepsilon_{a_i})$ by $I \cap D(\varepsilon_{a_i})$. Note that $A_i$ does not depend on the choice of $a_i$.
\end{defn}

We define dually the categories $B_j$.

\textbf{Example}:
Let Q be the quiver
$$
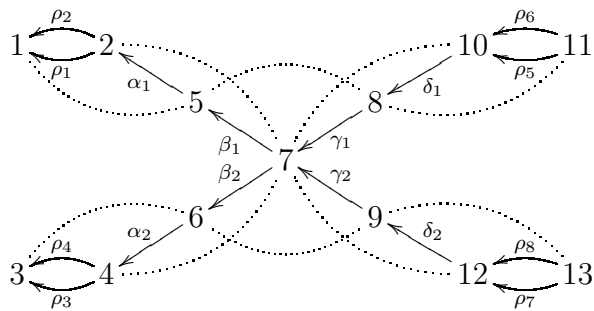
\xymatrix@R=10pt@C=10pt{
1 && 2 \ar@/^/[ll]^{\rho_1}\ar@/_/[ll]_{\rho_2}&&&&&& & & 10\ar[lld]^{\delta_1}\ar@{.}@/_0.5cm/[lllldd]&& 11\ar@/^/[ll]^{\rho_5}\ar@/_/[ll]_{\rho_6}\ar@{.}@/^0.5cm/[lllld]\\
 &&&&5\ar@{.}@/^0.5cm/[llllu]\ar[llu]^{\alpha_1}&&&&8\ar[lld]^{\gamma_1}\ar@{.}@/_0.5cm/[llll]&&&&  \\
 &&&&&&7\ar@{.}@/_0.5cm/[lllluu]\ar[llu]^{\beta_1}\ar[lld]_{\beta_2}\ar@{.}@/^0.5cm/[lllldd]&&&&&&  \\
 &&&&6\ar[lld]_{\alpha_2}\ar@{.}@/_0.5cm/[lllld]&&&&9\ar[llu]_{\gamma_2}\ar@{.}@/^0.5cm/[llll]&&&&  \\
3&& 4 \ar@/^/[ll]^{\rho_3}\ar@/_/[ll]_{\rho_4} && && &&& & 12\ar[llu]_{\delta_2}\ar@{.}@/^0.5cm/[lllluu]&& 13 \ar@/^/[ll]^{\rho_7}\ar@/_/[ll]_{\rho_8}\ar@{.}@/_0.5cm/[llllu]\\
}$$
\\with $I$ the ideal generated by $\alpha_1\rho_1$,
$\alpha_2\rho_4$, $\rho_5\delta_1$, $\rho_8\delta_2$,
$\beta_i\alpha_i$ for $i$ such that $i\in\{1,2\}$,
$\gamma_i\beta_i$ for $i$ such that $\in\{1,2\}$ and
$\delta_j\gamma_j=0$ for $j$ such that $j\in\{1,2\}$. Then $R=Q/I$
contains no DOZE. We have $B=B_1\times B_2$ with $B_1$ the full subcategory generated by $\{1,2,5\}$ and $B_2$ the full subcategory generated by  $\{3,4,6\}$. We also have $A=A_1\times A_2$ with
$A_1$ the full subcategory generated by $\{8,10,11\}$ and $A_2$
the full subcategory generated by $\{9,12,13\}$.

\begin{lem}
The categories $D(\varepsilon_{a_i})$ are full in $R$.
\end{lem}
Proof: By contradiction, suppose that $x$ and $y$ are objects of $A_i$ and let $\gamma: x \rightarrow y$ be an arrow which is not in $D(\varepsilon_{a_i})$. Then there exist strings $\omega: a_i \rightsquigarrow x$ and $\omega': a_i \rightsquigarrow y$, and since $\gamma$ is not in $A_i$, $\omega$ and $\omega'$ are such that $\omega \gamma$ and $\omega'\gamma^{-1}$ are reduced walks containing a zero-relation. But in this case, the zero-relation must contain $\gamma$ and thus $\omega= u\alpha_1...\alpha_n$ and $\omega'= u'\alpha_{m}^{-1}...\alpha_{n+2}^{-1}$, with $\alpha_1...\alpha_n\gamma$ and $\gamma\alpha_{n+2}...\alpha_m$ some relations. Thus, $\gamma\omega'^{-1}\Theta_i\rho'_{(i,\alpha)}$ is a
DOZE, a contradiction. We can apply the same proof to every morphism in the category $R$.
\\$\square$

\begin{lem}\label{aiaGauche}
Let $x \in (Q_{A_i})_0$ and $y \notin (Q_{A_i})_0$. Then there is no arrow $\alpha:y \rightarrow x$. As a consequence, if a relation $\rho$ does not start in $(Q_{A_i})_0$ then it does not end in $(Q_{A_i})_0$.
\end{lem}
Proof: Suppose that such an arrow exists. Let $\omega$  be a string from $a_i$ to $x$. Then $\omega\alpha^{-1}$ is not a reduced walk or else contains a zero-relation. If it is not a reduced walk, then $\omega=\omega'\alpha$, and $y \in (Q_{A_i})_0$, a contradiction. If it contains a zero-relation, then $\omega=\omega'\beta_n^{-1}...\beta_1^{-1}$ with $\alpha\beta_1...\beta_n$ in $I$ and so $\alpha\beta_1...\beta_n\omega'^{-1}\rho'_{i,\gamma}$, with $\gamma$ an exiting arrow of $\Theta_i$, is a DOZE, another contradiction.
\\$\square$

\begin{cor} The categories $D(\varepsilon_{a_i})$ are convex in $R$.
\end{cor}
Proof: This follows directly from lemma \ref{aiaGauche}.
\\$\square$

\begin{lem}\label{aucunCycle}
Each of the categories $D(\varepsilon_{a_i})$ contains only one simple cycle. In particular, it has no oriented cycle. Therefore, no relation $\rho$ starts in $(Q_{A_i})_0$ and ends outside the cycle $\Theta_i$.
\end{lem}

Proof: Suppose that there are two distinct cycles $\Theta$ and $\Theta_i$ in $D(\varepsilon_{a_i})$. If $\Theta$ contains a relation, we have a DOZE by gluing the relation with a walk relating the two cycles and which is chosen such that the two relations point in the same direction (this is possible because $\Theta$ and $\Theta_i$ are cycles). If not, we have a cycle with an exiting or an entering arrow and so we have a relation of the form $\alpha^-\alpha$ or of the form $\beta\beta^+$. In each case, we can construct a DOZE. The last statement follows from Lemma \ref{aiaGauche}.
\\$\square$

Recall that an algebra $R$ is left glued if $\mathcal R_R$ is cofinite in ind$R$ (see \cite{AC94}).

\begin{prop}\label{AlgebreInclinee} For all $i$, with $1
\leq i \leq n$, the algebra $A_i$ is a tilted algebra containing
a complete slice in its postprojective component. In
particular, it is left glued and all but a finite number of isomorphism classes of indecomposable $A_i$-modules
are in $\mathcal R_{A_i}$.
\end{prop}
Proof: Let $\alpha$ be an arrow from $x$ to $y$ in $A_i$. We want to show that the projective $A_i$-modules $P_x$ and $P_y$, with respective tops $S_x$ and $S_y$, are in the same component of the Auslander-Reiten quiver of $A_i$. The radical of $P_x$ is the direct sum of at most two terms, one of them is a uniserial module $L$ with top $S_y$. Moreover, there is an epimorphism $f$ from $P_y$ to $L$. Since the inclusion of the radical of a projective module into this projective module is irreducible, we only have to show that $f$ lies in a finite power of the radical of the module category.
\\Suppose $f_1:M\rightarrow L$ and $f_2:P_y\rightarrow M$ is a factorisation of $f$, where $M$ may be decomposable.
$$\xymatrix@R=10pt@C=10pt{
P_y\ar[rr]^{f_2} && M\ar[rr]^{f_1} && L\ar[rr] && P_x}.$$
Since $f$ is an epimorphism, $L$ is a quotient of $M$. Let $z_1$ belong to the support of $L$, $z_2$ belong to the support of $M$ but not to that of $L$ and $\beta$ be an arrow from $z_1$ to $z_2$. Since $L$ is a direct summand of rad$P_x$, there exists a walk $\omega'$ such that $\omega'\beta$ is in the support of $M$ and $\alpha\omega'\beta \in I$ ($\omega'$ may be trivial).
\\We have a relation of $I$ starting at $x$, thus $x$ is on the cycle $\Theta_i$ by Lemma \ref{aucunCycle}. Since there is no relation on a band, $\beta$ is not on a cycle (by Lemma \ref{aucunCycle}), and every arrow on a reduced walk beginning with $\alpha\omega'\beta$ and being the successor of $\beta$ is not in $\Theta_i$.
\\Thus, every arrow of the support of $M$ lies in the support of $L$ or is not in $\Theta_i$. Since $\Theta_i$ is not included in the support of $L$, it is not in the support of $M$.
\\The $k$-dimensions of the direct summands of $M$ are bounded above since the support of $M$ contains no band.
\\We showed that there exist only a finite number of isomorphism classes of indecomposable modules which are direct summands of a module through which $f$ factors.
\\Now, we only have to see that $A_i$ contains no double-zero (otherwise we have a DOZE since we have a string between every point of $A_i$ and $\Theta_i$). So $A_i$ is tilted (see \cite{HL00}) and we have our result by Theorem 3.4 of \cite{HL00}.
\\$\square$

\begin{defn}
Let $E$ be the set of strings which are not in one of the subcategories $A_i$ or $B_j$. We define the middle part $C$ to be the subcategory which has as set of objects $\bigcup_{\omega \in E}(D(\omega))_0$ and as set of morphisms $\Sigma_{\omega \in E}(D(\omega)(x,y))$ from $x$ to $y$.
\end{defn}

\begin{rem}
\label{remarqueCfini}
The algebra $C$ is of finite representation type since it contains no band and it is a string algebra (see \cite{BR87}).
\end{rem}

\subsection{Module Category}

We now begin the study of the category of indecomposable modules of $R$, with the aim to prove that it is a laura algebra under assumptions of this subsection. We denote by $A=\prod_{i=1}^n A_i$ the emph{right side algebra} of $R$ and by $B=\prod_{j=1}^m B_j$ the emph{left side algebra} of $R$.

\begin{prop}
$\emph{ind}R=\emph{ind}A\cup\emph{ind}B\cup\emph{ind}C$
\end{prop}
Proof: Every $A_i$-module or $B_j$-module can be considered as an $R$-module since $A_i$ and $B_j$ are full and convex subcategories of $R$, by completing their representations by zeros. Every $C$-module is a string module and thus is a string $R$-module.
\\Let $M$ be an indecomposable $R$-module and $\omega$ be a string lying in the support of $M$. Then each point and each arrow of the support of $M$ are in $D(\omega)$. Suppose that all points of $M$ are in one of the $A_i$, or in one of the $B_j$, but not all in the same. Let $x$ be a point of $A_i$ which is not in $B_j$ and $y$ be a point of $B_j$ which is not in $A_i$. Then the string $\omega: x \rightsquigarrow y$ cannot be in $A_i$, nor in $B_j$, and thus by definition lies in $C$.
$\square$

\begin{prop}\label{fermeSucc}
Let $M$ be an indecomposable $A_i$-module, $N$ be an indecomposable $R$-module
which is not an $A_i$-module and $f:M \rightarrow N$ be a non-zero morphism. Then
$M$ is a $C$-module or $M$ is a $B$-module.
\end{prop}
Proof: Suppose that $M$ is neither a $B$-module, nor a $C$-module. Then the support of $M$ contains at least a point $z'$ which is not an object of $B\cup C$. So it is an object of $A_i$ for some $i$. \\The image of $f$ is a submodule of $N$, so there exist a string $\omega$ going from $z'$ to $x$ and an arrow $\alpha:y\rightarrow x$, with $x$ in the support of the image of $f$ and $y$ in the support of $N$ but not in the support of the image of $f$.
\\Since $N$ is not an $A$-module, the support of $N$ contains at least a point $z$ which is not an object of $A_i$. Since $z$ and $y$ are in the support of $N$, there exists a string $\nu$ going from $y$ to $z$ of the form.
$$\xymatrix@R=10pt@C=10pt{
z'\ar@{~}[rrrr]^{\omega}&&&&x&&y\ar[ll]^{\alpha}\ar@{~}[rrrr]^{\nu}&&&&z}.$$
But $z'$ is not in $B\cup C$ and $z$ is not in $A_i$. Then there is no string between them, since such a string would be neither in $A$, nor in $B$, so it would be in $C$. If $\omega\alpha^{-1}\nu$ is a reduced walk, it contains a relation and this relation must contain $\alpha$. But in this case, we obtain a relation entering in $A_i$, a contradiction.
\\Thus $\omega\alpha^{-1}\nu$ is not reduced and the relation between $z$ and $z'$ goes in the other direction. Suppose that $\alpha_1\alpha_2...\alpha_n$ is going from $z'^*$ to $z^*$, as illustrated bellow.
$$\xymatrix@R=10pt@C=10pt{
z'^*\ar[rr]^{\alpha_1}\ar@{.}@/^0.5cm/[rrrrrrrr]&&...\ar[rr]^{\alpha_i}&&\cdot\ar@{~}[d]^{\eta}\ar[rr]^{\alpha_{i+1}}&&...\ar[rr]^{\alpha_n}&&z^*\\
&&&&x&&&&\\
&&&&y\ar[u]&&&&}.$$
Since $\omega$ and $\nu$ are strings, so are $\alpha_1...\alpha_i\eta$ and $\eta^{-1}\alpha_{i+1}...\alpha_n$. Since $R$ is a string algebra, this implies that $n=2$.
$$\xymatrix@R=10pt@C=10pt{
z'^*\ar[rr]^{\alpha_1}\ar@{.}@/^0.5cm/[rrrr]&&\cdot\ar@{~}[d]^{\eta}\ar[rr]^{\alpha_{2}}&&z^*\\
&&x&&\\
&&y\ar[u]&&}.$$
The target of $\alpha_1$ is on $\omega$ and on $\nu$, so it is in the support of Im$f$. Since Im$f$ is a quotient of $M$, $z'^*$ is in Im$f$. And since Im$f$ is a submodule of $N$, $z^*$ is in Im$f$, a contradiction to the fact that $\alpha_1\alpha_2 \in I$ and to Lemma \ref{aucunCycle}.
\\$\square$

\textbf{Proof of Theorem \ref{TheoA}:}
The first part of Theorem \ref{TheoA} follows now directly from \cite{HL00} (Theorems 2.6 and 3.4) and \cite{H01} (Theorem 3.1).

For the second part, the last proposition with the fact that
\\$|$ind$A_i\cap$ ind$B_j|\leq |$ ind$(A_i\cap B_j)|<\infty$ say that there exists only a finite number of isomorphism classes of indecomposable $A_i$-modules which admit successors from ind$R \setminus$ ind$A_i$. Let $\mathcal X_i$ be the subcategory of ind$A_i$ which contains modules without successors from ind$R \setminus$ ind$A_i$, that is
$$\mathcal X_i=\{X\in \textit{ind}A_i\textit{ | for all }0\neq f:X\rightarrow Y\textit{ in ind}R\textit{, }Y \in \textit{ ind}A_i\}.$$ We have that $\mathcal X_i$ is such that if $M$ is an $A_i$-module in $\mathcal X_i$ and if there is a non-zero morphism $f$ from $M$ to $M'$, then $M'$ is in $\mathcal X_i$. We define dually $\mathcal Y_j$.

We have that if $f:M \rightarrow N$ is a morphism, where $M$ and $N$ are objects of $\mathcal X_i$, then $f$ is irreducible in $A_i$ if and only if $f$ is irreducible in $R$, since the morphisms are preserved (this follows from the fullness of the categories $A_i$ and $B_j$ and from the fact that there are no points lying in $A_i \cap A_j$ if $i \neq j$). By studying the Auslander-Reiten translation in $A_i$ and in $R$, we can show that it is preserved. Take $$\xymatrix@R=10pt@C=10pt{
(*): &&0\ar[rr]&&M\ar[rr]&&_RI_0\ar[rr]^{i_1}&&_RI_1}$$
a minimal injective resolution of $M$ in $R$. We know that
$\tau_R^{-1}(M)=\mbox{Coker }\nu^{-1}(i_1)$, where $\nu$ denotes the Nakayama functor.
$$\xymatrix@R=10pt@C=10pt{
\nu^{-1}(_RI_0)\ar[rr]^{\nu^{-1}(i_1)}&&\nu^{-1}(_RI_1)\ar[rr]&&\mbox{Coker }\nu^{-1}(i_1)}$$
We compute $\tau_{A_i}^{-1}(M)$ in the same way.
In general, $\nu^{-1}(_{A_i}I_0) $ is a quotient of $\nu^{-1}(_RI_0)$ and $\nu^{-1}(_{A_i}I_1)$ is a quotient of $\nu^{-1}(_RI_1)$. The kernels $K_0$ and $K_1$ of those projections contain no points of $A_i$ in their supports.
$$\xymatrix@R=15pt@C=15pt{
0\ar[d]&&0\ar[d]&&\\
K_0\ar[d]\ar[rr]&&K_1\ar[d]\ar[rr]&&L\ar[d]\ar[rr]&&0\\
\nu^{-1}(_RI_0)\ar[rr]^{\nu^{-1}(i_1)}\ar[d]&&\nu^{-1}(_RI_1)\ar[rr]\ar[d]&&\tau^{-1}_R(M)\ar[d]^g\ar[rr]&&0\\
\nu^{-1}(_{A_i}I_0)\ar[rr]^{\nu^{-1}(i_1)}\ar[d]&&\nu^{-1}(_{A_i}I_1)\ar[rr]\ar[d]&&\tau^{-1}_{A_i}(M)\ar[rr]&&0\\
0&&0&&}$$
Let $L$ be the cokernel of the induced morphism from $K_0$ to $K_1$. It contains no points of $A_i$. But $\tau^{-1}_R(M)$ is an $A_i$-module and so $L$ is zero and $g$ is a monomorphism. It is also an epimorphism.
Thus a path in $\mathcal X_i$ is sectionally refinable in $\mathcal X_i$ if and only if it is also in ind$R$.

Moreover, let $P$ be an $A_i$-modules in $\mathcal X_i$, if $P$ is $R$-projective, then it is $A_i$-projective. If $P$ is a $R$-projective module, than for every $R$-epimorphism, and in particular for every $A_i$-epimorphism $f:M\rightarrow N$ with $g: P\rightarrow N$, there exists an $R$-morphism $h : P\rightarrow M$ such that $hf=g$. Since $h$ is a morphism between two $A_i$-modules, it is a morphism in mod$A_i$. On the other hand, an $A_i$-module $I$ in $\mathcal X_i$ is an injective $R$--module if and only if $I$ is an injective $A_i$-module (this follows from Lemma \ref{aiaGauche} and from the construction of injective modules).

Now, let $(Q,I)$ be a string bound quiver having a band and such that each band has only exiting arrows or only entering arrows. Let $R=kQ/I$ have no DOZE. We show that if $M$ is not in $\mathcal L_R \cup \mathcal R_R$, then one of the following conditions is satisfied:
\begin{enumerate}
\item The module $M$ is in ind$C$;
\item There exist $i$ and $j$ such that $1 \leq i \leq n$ and $1 \leq j \leq m$, and such that the module $M$ is in ind$A_i \cap$ind$B_j$;
\item The module $M$ is in ind$A_i$, but not in $\mathcal R_{A_i}$;
\item The module $M$ is in ind$B_j$, but not in $\mathcal L_{B_j}$.
\end{enumerate}
If $M$ does not satisfy the first two conditions, then by Proposition \ref{fermeSucc} there exists $i$ such that $M$ is in $\mathcal X_i$ or $j$ such that $M$ is in $\mathcal Y_j$. In the first case, if $M$ is not in $\mathcal R_R$, then there exists a path of $R$-morphisms which is not sectionally refinable from $M$ to an $R$-projective indecomposable module. From the results above, we obtain that $M$ is not in $\mathcal R_{A_i}$. In the second, if $M$ is not in $\mathcal L_R$, then there exists a path of $R$-morphisms which is not sectionally refinable from an $R$-injective indecomposable module to $M$. We obtain that $M$ is not in $\mathcal L_{B_j}$. We have shown our statement.

In the four cases, $M$ belongs to a finite set and thus all but a finite number of isomorphism classes of modules of ind$R$ are in $\mathcal L_R \cup \mathcal R_R$. In particular, $R$ is laura.
It cannot be quasi-tilted of canonical type since none of its bands have entering and exiting arrows (see Theorems 2.6 and 3.4 of \cite{HL00}).
\\$\square$

\begin{cor}\label{Domestique} Let $(Q, I)$ be a string bound quiver having bands and such that each band has only exiting arrows or only entering arrows. Then, if $R=kQ/I$ is a strict laura or a tilted string algebra without DOZE, its left and right end algebras are tilted of type $\widetilde{\mathbb{A}}_n$. Consequently, each string strict laura or tilted algebra without DOZE is domestic.
\end{cor}

Proof: We have that $\mathcal X_i$ is such that if $M$ is in $\mathcal X_i$ and there exist a non-zero morphism from $M$ to $M'$, then $M'$ is in $\mathcal X_i$, and such that only a finite number of $A_i$-modules are not in $\mathcal X_i$. Therefore, if $R$ is not quasi-tilted, we have that $\mathcal X_i$ contains a complete slice of mod$A_i$ and this complete slice gives us a connected component of the right end algebra of $R$. Moreover, the right end algebra cannot contain another factor (in this case, by the arguments of the proof of Theorem \ref{TheoA}, we obtain an infinite number of $R$-modules which are not in $\mathcal L_R \cup \mathcal R_R$). The domesticity of $R$ follows from \cite{AC03}.
\\$\square$

\section[DOZED string modules]{DOZED string modules}
We now study string algebras whose bound quiver contains at least a DOZE.

\begin{defn}
Let $R=kQ/I$ be a string algebra with a DOZE
$$\alpha_1...\alpha_l\omega_1{\omega_2}^n\omega_3\beta_m...\beta_1$$
for some walks $\alpha_1\ldots\alpha_l$, $\omega_1$, $\omega_3$, $\beta_m...\beta_1$ and a band $\omega_2$.
The string module $M_n = M(\sigma_n)$ corresponding to the string
$$\sigma_n=\alpha_3...\alpha_l\omega_1{\omega_2}^n\omega_3\beta_m...\beta_{3}$$
is called the DOZED module of power $n$. By convention, $\alpha_3...\alpha_l$ represents the trivial path if the relation $\alpha_1...\alpha_l$ is of length 2.
\end{defn}

\begin{thm}\label{TheoC}
Let $R$ be a string algebra containing a DOZE. Then the DOZED string modules have projective and injective dimension greater than one. Thus, $R$ is not laura.
\end{thm}
Proof: Let $\rho_1\omega_1\omega_2\omega_3\rho_2$ be a DOZE where $\omega_2$ is a non-trivial band, $\rho_1=\alpha_1...\alpha_l$ and $\rho_2=\beta_m...\beta_1$. Let $M_n=M(\sigma_n)$ be the DOZED module and $I$ the injective envelope of $M_n$.
\\Let $\alpha_3...\alpha_l\alpha_{l+1}...\alpha_k$ be the longest path such that $\sigma_n=\alpha_3...\alpha_l\alpha_{l+1}...\alpha_k\sigma'_n$ and $x_3$, $x_4$, ... $x_{k+1}$ the vertices of this path. Then $x_{k+1}$ is a sink of $\sigma_n$ and thus $I_{x_{k+1}}$ is a direct summand of $I$. Let $\pi_1$ be the canonical projection on $I_{x_{k+1}}$, $M_n(x)$ and $M_n(\alpha)$ be the vector spaces and the linear maps corresponding to the vertex $x$ and the arrow $\alpha$ in the representation $M_n$, respectively and $f_{x_i}$ be the linear map associated to $x_i$ in the representation of the morphism $f$. Then we have the following identities: \begin{itemize}
\item $\pi_1\circ M_n(\alpha_i)={\pi_1}_{x_i}$ for every $i$ such that $3 \leq i \leq k$;
\item $\pi_1\circ M_n(\alpha_2)=0$ by definition of $M_n$;
\item $\pi_1\circ I(\alpha_i)={\pi_1}_{x_i}$ for every $i$ such that $2 \leq i \leq k$, since $\alpha_i...\alpha_k$ is the first non-zero path going to $t(\alpha_k)$ which is a sink of $\sigma_n$;
\item $\pi_1\circ I(\alpha_1)=0$ since $\alpha_1...\alpha_l$ is a relation ending in $x_{l+1}$ and of minimal length for this property.
\end{itemize}
Let $\iota$ be the inclusion of $M_n$ in $I$, $\iota_x$ the linear map induced by $\iota$ between $M_n(x)$ and $I(x)$, $C=$Coker$ \iota$ and $c=$coker$ \iota$. Then:
\begin{itemize}
\item $\pi_1\circ \iota_{x_i}={\pi_1}_{x_i}$ for every $i$ such that $3 \leq i \leq k$;
\item $\pi_1\circ \iota_{x_2}=0$, by applying $\pi_1$ to the equation $\iota_{x_3}\circ M_n(\alpha_2)=I(\alpha_2)\circ \iota_{x_2}$; \item
$\pi_1(C(x_2))\cong k$ since $\pi_1\circ \iota_{x_2}=0$ and since
$\pi_1(I(x_2))\cong k$;
\item $\pi_1 \circ c_{x_2}={\pi_1}_{x_i}$;
\item $\pi_1 \circ c_{x_3}=0$;
\end{itemize}

$$\xymatrix@R=10pt@C=30pt{
&\cdot\ar[rr]\ar[ddd]&&\cdot\ar[rr]^{0}\ar[ddd]^{0}&&
\cdot\ar[rr]^{id_k}\ar[ddd]^{id_k}&&\cdot\ar[ddd]^{id_k}\\
M_n(x_1)\ar[rr]^{M_n(\alpha_1)}\ar[ddd]^{\iota_{x_1}}\ar[ru]^{\pi_1}&&
M_n(x_2)\ar[rr]^{M_n(\alpha_2)}\ar[ddd]^{\iota_{x_2}}\ar[ru]^{\pi_1}&&
M_n(x_3)\ar[rr]^{M_n(\alpha_3)}\ar[ddd]^{\iota_{x_3}}\ar[ru]^{\pi_1}&&
M_n(x_4)\ar[ddd]^{\iota_{x_4}}\ar[ru]^{\pi_1}&\\
&&&&&&&\\
&\cdot\ar[rr]^{0 \qquad}\ar[ddd]&&\cdot\ar[rr]^{id_k
\qquad}\ar[ddd]^{id_k}&&
\cdot\ar[rr]^{id_k \qquad}\ar[ddd]^{0}&&\cdot\ar[ddd]\\
I(x_1)\ar[rr]^{I(\alpha_1)}\ar[ddd]^{c_{x_1}}\ar[ru]^{\pi_1}&&
I(x_2)\ar[rr]^{I(\alpha_2)}\ar[ddd]^{c_{x_2}}\ar[ru]^{\pi_1}&&
I(x_3)\ar[rr]^{I(\alpha_3)}\ar[ddd]^{c_{x_3}}\ar[ru]^{\pi_1}&&
I(x_4)\ar[ddd]^{c_{x_4}}\ar[ru]^{\pi_1}&\\
&&&&&&&\\
&\cdot\ar[rr]^{0 \qquad}&&k \ar[rr]^{0 \qquad}&&\cdot\ar[rr]&&\cdot\\
C(x_1)\ar[rr]^{C(\alpha_1)}\ar[ru]^{\pi_1}&&C(x_2)\ar[rr]^{C(\alpha_2)}\ar[ru]^{\pi_1}&&
C(x_3)\ar[rr]^{C(\alpha_3)}\ar[ru]^{\pi_1}&&C(x_4)\ar[ru]^{\pi_1}&\\
}$$

By applying $\pi_1$ to the equation $C(\alpha_1)\circ c_{x_1}=c_{x_2}\circ I(\alpha_1)$, we obtain that $\pi_1\circ
C(\alpha_1)\circ c_{x_1}=\pi_1 \circ c_{x_2} \circ I(\alpha_1)=\pi_1 \circ I(\alpha_2)=0$. Since $c_{x_1}$ is an epimorphism, we have $\pi_1 \circ C(\alpha_1)=0$. We show in the same way that $\pi_1 \circ C(\alpha_2)=0$, and that $\pi_1 \circ C(\beta)=0$ for every arrow $\beta$ having $x_2$ for its source.
\\Thus $C$ admits an indecomposable direct summand whose support does not contain $x_1$ and admits $x_2$ as a sink. So this direct summand (and thus $C$) is not an injective module.
\\$\square$

\textbf{Proof of Theorem \ref{TheoPresqueFinal}:}
a) implies b):
\begin{enumerate}
\item If there is no band, the algebra is of finite representation type by \cite{BR87}, and we are done;
\item Otherwise, Theorem \ref{TheoA} gives us the statement.
\end{enumerate}
The statement b) implies trivially c). The statement c) implies a)
by Theorem \ref{TheoC}.
\\$\square$

We now generalize Theorem \ref{TheoPresqueFinal} to special biserial algebras.

\section{Laura special biserial algebras and Skowro{\'n}ski's conjecture}

\begin{defn}
\label{biserielleSpeciale}
Let $R$ be an algebra. It is a special biserial algebra if it admits a presentation $R=kQ/I$ such that : \begin{enumerate}
\item Each point has at most two arrows entering and two arrows exiting;
\item For each arrow  $\alpha:x \rightarrow y$ there is at most one arrow $\beta:y \rightarrow z $ such that $\alpha\beta$ is not in $I$ and at most one arrow $\gamma: z \rightarrow x$ such that $\gamma\alpha$ is not in $I$.
\end{enumerate}
\end{defn}

\underline{Remark:} For a special biserial algebra $R$, we denote $J$ the ideal generated by the paths appearing in the commutativity relations of $I$. Then, $R/J$ is a string algebra.

\begin{defn}
Let $R=kQ/I$ be a special biserial algebra and $(\rho_1,\rho_2)$ be two consecutive zero-relations on $Q$. We say that they form a DOZE if they do so on $R/J$.
\end{defn}

By \cite{WW85}, if $R$ is special biserial, any indecomposable $R$-module which is not projective-injective is also an $R/J$-module. This implies that any indecomposable $R$-module which is not projective-injective is a string module or a band module. Moreover, the restriction of scalars functor mod$R/J \rightarrow$ mod$R$ is full, faithful and sends an irreducible morphism in mod$R/J$ onto an irreducible morphism in mod$R$.

\textbf{Proof of Theorem \ref{TheoFinal}:}
 The proof that (b) implies (c) follows from the definition of laura algebras.
\\Suppose that there exists a DOZE on $(Q,I)$ and let $M_n$ be a DOZED module on $R/J$. Then $M_n$ is also an $R$-module. We have shown in the proof of Theorem \ref{TheoC} that the injective dimension over $R/J$ of $M_n$ is greater than or equal to two by building $I_{x_{l+1}}$ and the cokernel $C$. Remark that a DOZE on $R$ contains no path involved in a binomial relation. By applying the same technique, we show that the injective dimension of $M_n$ over $R$ is greater than or equal to two and we get that c) implies a).
\\The proof that a) implies d) follows from Theorem \ref{TheoC}.
\\For the last implication, let $R$ be a special biserial algebra such that $R/J$ is laura. Let $M$ be a non
projective-injective $R$-module which is not in $\mathcal L_R \cup \mathcal R_R$. Then there exist a non sectionally refinable path (*) going from $M$ to a $R$-projective module $P$ and a non sectionally refinable path (**) going from an $R$-injective module $I$ to $M$.
$$\xymatrix@R=10pt@C=10pt{
(*)&M\ar[rd] &&\tau^{-1} M \ar@{~>}[rr]&& P &&(**)&I\ar@{~>}[rr] &&\tau M\ar[rd] &&M \\
&&\cdot \ar[ru]&&&&&&&&&\cdot \ar[ru]&}$$
\\The fact that (*) and (**) are not sectionally refinable is preserved in $R/J$ since an almost split sequence admitting a projective-injective middle term has at least another middle term \cite{ASS97}. Moreover, if $P$ and $I$ are not projective-injective modules, then they are $R/J$-projective and $R/J$-injective modules respectively. If $P$ is projective-injective, then $P/$soc$P$ is an $R/J$-projective module. We only have to compose (*) with the left minimal almost split morphism going from $P$ to $P/$soc$P$, and we get the desired path from $M$ to an $R/J$-projective module.
$$\xymatrix@R=10pt@C=10pt{
M\ar@{~>}[rr]&& P\ar[r]&P/\emph{soc}P }$$
In all cases, if $M$ is a non projective-injective $R$-module which is not in $\mathcal L_R \cup \mathcal R_R$, then it is not in $\mathcal L_{R/J} \cup \mathcal R_{R/J}$, and so the set of indecomposable $R$-modules which are not in $\mathcal L_R \cup \mathcal R_R$ is finite.
\\$\square$

We conclude with an example.

\textbf{Example:} Let $R=kQ/I$, where $Q$ is the following quiver:
$$
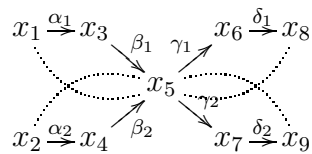
\xymatrix@R=10pt@C=10pt{
x_1\ar@{.}@/_0.5cm/[rrd]\ar[r]^{\alpha_1}& x_3\ar[rd]^{\beta_1}&& x_6\ar[r]^{\delta_1} &x_8  \\
&&x_5\ar[ru]^{\gamma_1}\ar[rd]^{\gamma_2}\ar@{.}@/^0.5cm/[rrd]\ar@{.}@/_0.5cm/[rru]& &&\\
x_2\ar@{.}@/^0.5cm/[rru]\ar[r]^{\alpha_2}& x_4\ar[ru]_{\beta_2}&& x_7\ar[r]^{\delta_2}&x_9
}$$
\\and $I$ is generated by $\alpha_i\beta_i$ and $\gamma_i\delta_i$.
This algebra is not a string algebra, but by applying the evident action of the group $\mathbb Z_2$, we obtain the skew group algebra $A[G]$ given by the following quiver:
$$
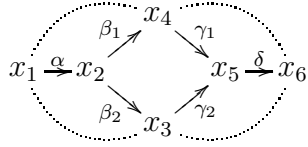
\xymatrix@R=10pt@C=10pt{
            &                                          & x_4\ar@{.}@/^0.5cm/[rrd]\ar[rd]^{\gamma_1} & &   \\
x_1\ar@{.}@/^0.5cm/[rru]\ar@{.}@/_0.5cm/[rrd]\ar[r]^{\alpha}   & x_2\ar[ru]^{\beta_1}\ar[rd]_{\beta_2} &                        &  x_5\ar[r]^{\delta} & x_6\\
  &                                        & x_3\ar@{.}@/_0.5cm/[rru]\ar[ru]_{\gamma_2}&&
}$$
\\where $I$ is generated by $\alpha\beta_1$, $\alpha\beta_2$, $\gamma_2\delta$ and $\gamma_1\delta$. It is a string algebra which contains the DOZE $\alpha\beta_1\gamma_1\delta$, this latter algebra is not laura by Theorem \ref{TheoPresqueFinal} and so the initial algebra $R$ is not laura by \cite{ALR07}.

%-------------------------------------------------------------------------------
%
%   ACKNOWLEDGEMENTS
%
%-------------------------------------------------------------------------------
%\bigskip
\noindent {\bf ACKNOWLEDGEMENTS.} The author warmly thanks Ibrahim Assem as well as Patrick Le Meur and Shiping Liu for their suggestions. She also thanks ISM and Sherbrooke University for their financial support.

%-------------------------------------------------------------------------------
%
%  BIBLIOGRAPHIE
%
%-------------------------------------------------------------------------------

\bibliography{bibliolaura}

\end{document}